\documentclass{ifacconf}

\usepackage{graphicx}      
\usepackage{natbib}        

\usepackage{amsmath, amssymb}
\usepackage{algorithm}
\usepackage{algpseudocode}
\usepackage{multirow}

\def\ag#1{{\color{black}#1}} 
\def\pars#1{{\color{black}#1}} 

\usepackage{xcolor}
\everymath{\displaystyle}

\begin{document}
\begin{frontmatter}

\title{ADMM-based Distributed State Estimation for Power Systems: Evaluation of Performance\thanksref{footnoteinfo}} 

\thanks[footnoteinfo]{The research of A. Gasnikov is partially supported by RFBR project 19-31-51001 and by the Ministry of Science and Higher Education of the Russian Federation (Goszadaniye) No. 075-00337-20-03, project No. 0714-2020-0005.}


\author[Skoltech]{Sergei Parsegov}
\author[MIPT]{Samal Kubentayeva} 
\author[Skoltech]{Elena Gryazina}
\author[MIPT,IITP,CMC]{Alexander Gasnikov}
\author[Skoltech]{Federico Ib\'a\~nez}

\address[Skoltech]{Skolkovo Institute of Science and Technology, Moscow, Russia (e-mail: e.gryazina@skoltech.ru, s.parsegov@skoltech.ru)}
\address[MIPT]{Moscow Institute of Physics and Technology, Moscow, Russia (e-mail: samal.kubentayeva@phystech.edu)}
\address[IITP]{A.A. Kharkevich Institute for Information Transmission Problems, Russian Academy of Sciences, Moscow, Russia}
\address[CMC]{Caucasus Mathematical Center, Adyghe State University, Maykop, Russia (e-mail: gasnikov.av@mipt.ru).}

\begin{abstract}                
Recently, distributed algorithms for power system state estimation have attracted significant attention. Along with such advantages as decomposition, parallelization of the original problem and absence of a central computation unit, distributed state estimation may also serve for local information privacy reasons since the only information to be transferred is the boundary states of neighboring areas. In this paper, we propose some novel approaches for speeding up the ADMM-based distributed state estimation algorithms by utilizing some recent results in optimization theory. We also thoroughly analyze the theoretical and practical performance, concluding that accelerated approach outperforms the existing ones. The theoretical considerations are verified through the experiments on a scalable example.
\end{abstract}

\end{frontmatter}

\section{Introduction}

State estimation (SE) is a well-known tool used to capture the real-time operating conditions of a power system (mainly bus voltage magnitudes and phase angles). Although the problem has a long history with the first results published in early 1970's, see \cite{Schweppe}, power system state estimation (PSSE) is now one of the most vital components providing proper data for operation and management of power systems, e.g. see \citet{bookGomez} and \citet{bookConejo}. Traditionally, each regional control center performs SE independently with very limited interaction between other control centers. However, novel solutions alternative to classical centralized frameworks are needed due, see e.g. \cite{kekatos2013distributed} and \cite{ZhengWu_et_al2017}:
\begin{itemize}
    \item the scale and complexity of networks are constantly growing. This increases the computational burden and produces communication bottlenecks;
    \item independent system operators are no longer fully independent due to more intensive use of the so-called tie-lines interconnecting the areas, that need to be monitored;
    \item the power transmission between areas becomes less predictable due to uncertain renewable energy sources that requires faster state estimation;
    \item the information privacy of subsystems should be preserved, therefore the areas cannot share their internal states: a centralized solution to SE is thus infeasible. Moreover, there are no coordinators above independent system operators.
\end{itemize}

For the problems mentioned above, the use of distributed optimization algorithms for PSSE continues to gain popularity. Advantages of distributed approaches include, but not limited to increased robustness as well as reduction in computation, communication, and memory requirements per area since each local control center requires only a subset of the global information. Having local computation units within each area, the implementation of distributed state estimation can replace a centralized coordinator: instead of collecting, storing and processing of huge measurement data for the state estimation, the same goal can be achieved by dealing with local measurements and their processing with the exchange of boundary state variables only. Besides preventing the appearance of communication bottlenecks, this may also look fascinating due to preserving information privacy if the interacting areas do not want to indicate their internal states. 

The most promising in terms of lower computation burdens, simplicity and better converge rate, is a distributed state estimation method proposed in \citet{kekatos2013distributed}, \citet{Kekatos}. Their approach is based on the alternating direction method of multipliers (ADMM). ADMM-based technique applied to secure operation and distributed state estimation method under hybrid cyber attacks was recently studied by \citet{Du}.  Some other problem formulations and corresponding distributed approaches utilizing ADMM can be found in \cite{Bolognani2014}, \cite{ZhengWu_et_al2017}.

Due to the constantly increasing complexity of networks and higher penetration of stochastic renewable generation, it is needed to perform distributed state estimation more frequently. The acceleration of existing distributed methods or/and the development of novel ones, along with increasing their robustness are of great practical importance, i.e. for blackouts prevention. Thus, robust and rapid distributed state estimator can be useful in extreme situations, so that critical parts of the network could be detected and proper actions taken to prevent rolling outages. Nevertheless, as compared to the existing results, the problem of computational speedup in distributed ADMM-based PSSE setup is still lacking.

During the past five years, a number of accelerated distributed methods have been proposed, e.g. see \cite{Lan}, \cite{lin2020accelerated} and \cite{dvinskikh2019decentralized}. Therefore, it is of great interest to test their applicability to state estimation problems and find out whether/when their implementation can provide some benefits.

In this paper, we perform rigorous analysis and comparison of ADMM-based distributed optimization approaches applied to a specific PSSE problem. The major contributions of this paper are summarized as follows: we
\begin{itemize}
    \item demonstrate that conventional ADMM for distributed PSSE can be slightly improved, guaranteeing faster convergence,
    \item adopt accelerated version of the method to the distributed case, study its rate of convergence and compare the performance of both approaches,
    \item analyze other aspects and possibilities of speeding up the convergence under additional assumptions on the network.
\end{itemize}

Our investigations show, that within the specific problem formulation we consider, with no other assumptions on network topology available, the accelerated ADMM-based approach is the best choice. Numerical simulations on synthetic scalable network topology with different numbers of interconnected areas confirm our theoretical discussions.

The remainder of the paper is organized as follows. In Section 2 some basic of SE are given and the problem formulation is stated. In Section 3 the distributed versions of ADMM and the accelerated ADMM (A-ADMM) are studied and compared, the corresponding convergence theorem is formulated. Then the applicability to distributed SE of other accelerated approaches is discussed. In Section 4 the examples are considered, that illustrate the performance of the methods and confirm the theoretical part. Finally, the conclusion sums up the objectives of the paper in addition to the major results. 

\section{Problem Formulation}

Power system state estimation (PSSE) is an important tool aimed to find unknown values of the state variables using some imperfect measurements. The classical nonlinear model for power systems is a set of over-determined equations, e.g. see \cite{AburBook}:
\begin{eqnarray}
    \label{eq:nonlin_model}
    \mathbf{z} = h(\mathbf{v}) + \boldsymbol{\varepsilon},
\end{eqnarray}
where $\mathbf{v}$ denotes the state vector, $\mathbf{z}$ is the vector of measurements, $h(\mathbf{v})$ is the vector of functions relating measurements to states, and $\boldsymbol{\varepsilon}$ is a disturbance term that captures errors and inaccuracies. For simplicity, the errors in vector $\boldsymbol{\varepsilon}$ are typically assumed zero mean with known diagonal covariance matrix. This matrix can be considered as an identity matrix after simple coordinate scaling.

The most commonly-used approach to estimate the actual value of the unknown variables is to solve the ``conventional'' least-squares optimization problem:
\begin{equation}
    \label{eq:optgeneral}
   \min_{\mathbf{v} } \frac{1}{2}\|\mathbf{z} - h(\mathbf{v})\|^2,
\end{equation}
which is usually non-convex. However, in some cases nonlinear model \eqref{eq:nonlin_model} reduces to a linear one: 
$$
   \mathbf{z} = H \mathbf{v} + \boldsymbol{\varepsilon},
$$ 
where measurement matrix $H$ is known.

Models of such kind appear, e.g., within DC approximation for power flow model, that is sufficiently accurate for high voltage modes in large-scale power grids. The other advantages of the model and its applicability are described in \citet{DC}.

Moreover, linear models for state estimation directly follow from the measurements provided by high-speed synchronized sensors: phasor measurement units (PMUs), e.g. see \cite{PMU2008}.


Similar to \cite{kekatos2013distributed}, we consider a power system that consists of $K$ interconnected control areas with their own local control centers. Each area can be considered
of as an independent system operator region, a balancing authority, a power distribution center, or a substation, for detail see \cite{Kekatos}. The local state vectors of neighboring areas $\mathbf{v}_k$, $k\in 1:K$ overlap partially. Without loss of generality we consider a linear model that corresponds to DC approximation approach. The model of each area is as follows:

\begin{eqnarray}\label{model_active_power}
 \mathbf{z}_k = 
 H_k \mathbf{v}_k + \boldsymbol{\varepsilon}_k, \quad \forall k \in K, 
\end{eqnarray}
where the vector $\mathbf{v}_k \in \mathbb{R}^{n_k}$ collects the system states related to area $k$,  $\mathbf{z}_k \in \mathbb{R}^{n_k + l_k}$ denotes the measurements from the sensors on buses and lines located in area $k$, $\boldsymbol{\varepsilon}$ is a zero mean Gaussian noise. The measurement matrix of the $k$-th area is denoted by $H_k \mathbf{v}_k$. Thus, we arrive at the following SE problem for the entire power system:
\begin{eqnarray}\label{optimization_prob}
    \min_{\mathbf{v}_k \in \mathcal{X}_k} f_k(\mathbf{v}_k) = \min_{\mathbf{v_k} \in \mathcal{X}_k} \frac{1}{2}\|\mathbf{z}_k - H_k \mathbf{v}_k\|^2,
\end{eqnarray}
where the convex set $\mathcal{X}_k$ collects prior information about the buses in the grid $\mathcal{G}$. For detailed exposition see in \citet{Monticelli}.

According to the structure of the network, each area has a common node (bus) with another one. Two areas that share a common bus restore the value of this bus, also they should take into account the value of each other and adjust shared state. Therefore, in order to solve the PSSE for the entire network, we need to introduce coupling conditions for the neighboring areas. Denote the voltage of a shared bus of areas $k$ and $l$ by $\mathbf{v}_k [l]$ and $\mathbf{v}_l [k]$. Thus, we obtain a consensus condition $\mathbf{v}_k [l] = \mathbf{v}_l [k] = \mathbf{u}_{kl}$, where $\mathbf{u}_{kl}$ is the new variable need for further construction of distributed methods. 

Finally, the optimization problem to be solved by distributed version of ADMM and A-ADMM looks as follows: 
\begin{eqnarray}\label{optim_prob_for_distr_methods}
&\min_{\{\mathbf{v}_k \in \mathcal{X}_k\}} \sum_{k=1}^{K} f_k(\mathbf{v}_k)& \notag\\
\hbox{s.t.}& \mathbf{v}_k[l]= \mathbf{u}_{kl} & \forall l \in \mathcal{B}_{k}, \notag\\
&\mathbf{v}_l[k]= \mathbf{u}_{kl} & \forall k \in \mathcal{B}_{l},
\end{eqnarray}
where $\mathcal{B}_k$ is the set of areas sharing states with area $k$, $\mathcal{B}_l$, respectively, with the same way with area $l$.

In the next section we analyze the conventional ADMM and its accelerated modification applied to distributed PSSE and evaluate their convergence rates.

\section{Methods: Distributed ADMM and A-ADMM}
In this section we focus on two efficient methods aimed to solve problem (\ref{optim_prob_for_distr_methods}) in a distributed fashion. Note that the methods are extensively used in power systems, see e.g. the works by \citet{kekatos2013decentralized} and \citet{LA-ADMM}.

The first method we discuss is the distributed ADMM derived and described in \citet{Boyd}. The second is the distributed modification of accelerated ADMM based on the results by \cite{lin2020accelerated}. We adopt both methods to the PSSE problem and estimate the convergence rate for the ergodic version of ADMM, that slightly outperforms the one proposed in \cite{kekatos2013distributed}, \cite{Kekatos} as applied to our problem, see Theorem~\ref{Theorem} for details. 
\subsection{Alternating Direction Method of Multipliers}
%
\begin{algorithm}
\caption{Alternating Direction Method of Multipliers}
\label{algADMM}
\begin{algorithmic}[1]
\Require  matrix $H$, available measurements $\mathbf{z}$, parameter $\mu$,tolerances $\varepsilon_{primal}=10^{-3}, \; \varepsilon_{dual}=10^{-4}$, initial point $\mathbf{u}_{kl}^0 = 0$, $\mathbf{y}_{kl}^0 = 0$ for each $k, l \in K$ 
\State $i \leftarrow 0$
\Repeat
\State \begin{eqnarray*}
       &\mathbf{v}_k^i &= \min_{\mathbf{v}_k \in K} f_k(\mathbf{v}_k) + \sum_{l \in B_k} \mathbf{y}_{kl}^{i-1}\left(\mathbf{v}_k[l] - \mathbf{u}^{i-1}_{kl}\right) \\
       &&+ \sum_{l \in B_k} \frac{\mu}{2} \|\mathbf{v}_k[l] - \mathbf{u}_{kl}^{i-1}\|_2^2
      \end{eqnarray*}
\State $\mathbf{u}_{kl}^i = \frac{1}{2}\left(\mathbf{v}_k^i[l] + \mathbf{v}_l^i[k]\right) $ 
\State $\mathbf{y}_{kl}^i = \mathbf{y}^{i-1}_{kl} + \mu \left( \mathbf{v}_k^i[l] - \mathbf{u}_{kl}^i\right)$ 
\Until {$\|\mathbf{r}_k^i\|_2^2 < \varepsilon_{primal}$ and $\|\mathbf{s}_k^i\|_2^2 < \varepsilon_{dual}$}
\State $i \leftarrow i+1$
\end{algorithmic}
\end{algorithm}

The major benefit of ADMM is that it may serve as a distributed solver, see Algorithm \ref{algADMM}. The objective function of the problem (\ref{optim_prob_for_distr_methods}) is separable and the variables $\mathbf{v}_k\in \mathbb{R}^{n_k}$ are subvectors of the system state vector; that is each of these local variables $\mathbf{v}_k$ consists of a selection of components of the global variable $\mathbf{v}$. Therefore, the $\mathbf{v}$-minimization step is split into $K$ separate problems that can be solved in parallel like the local area minimization problems in Algorithm \ref{algADMM}. In the $\mathbf{u}$-update step the neighboring areas collaborate (interact) each other to develop a global state vector. Hence, we divide by two because an element of the global state vector can be restored by two local  variables of the state regions that have their own opinion guess about it. In other words, we consider local state vectors that directly configure an element of the global vector.

Note that at this stage we are restoring elements that belong to the interacting areas.
Further, the  $\mathbf{y}_{kl}$-update step that can be interpreted as a vector of price adjustments for conditions of the problem (\ref{optim_prob_for_distr_methods}) and this step can be carried out independently in parallel for each $k\in K$. In addition, the dual variable update uses a step size equal to the augmented Lagrangian parameter $\mu$, where $\mu>0$.

The necessary and sufficient optimality conditions for the problem (\ref{optim_prob_for_distr_methods}) are primal and dual feasibility. Therefore, we check the termination criterion at step 6 where we introduce primal residual $\mathbf{r}_k^i = \mathbf{v}_k[l]^i - \mathbf{u}_{kl}^i$ for primal feasibility and dual residual $\mathbf{s}_k^i = \mathbf{u}_{kl}^i - \mathbf{u}_{kl}^{i-1}$ for dual feasibility. Note that when $\mathbf{r}^i$ and $\mathbf{s}^i$ are small, the objective suboptimality also must be small. Hence, a reasonable stopping criterion is that the dual and primal residuals must be small.

In general, ADMM solves convex optimization problems of the form
\begin{eqnarray}\label{Lan_gen_prob}
&& \min_{\mathbf{x}\in \mathcal{X}, \mathbf{y}\in \mathcal{Y}} f(\mathbf{x}) + g(\mathbf{y}) \notag \\
&& \hbox{s.t.} \; A\mathbf{x} + B \mathbf{y} = b, 
\end{eqnarray}
where matrices and vector of proper dimensions $(A, B, b)$ are given. See \citet{Lan} article for more details. The problem (\ref{optim_prob_for_distr_methods}) is a special case of the general problem (\ref{Lan_gen_prob}). 

One can notice, that the main difference between our problem (\ref{optim_prob_for_distr_methods}) and the general one (\ref{Lan_gen_prob}) is that we use the sum of functions of the same type with variables $\mathbf{v}_k \in X_k$ and in our case $g(\mathbf{y}) = 0$. Thus, we formulate a convergence theorem for our special case as applied to the PSSE problem, the proof of which uses a similar way as in Theorem 3.10 in \citet{Lan}. 

\begin{thm}\label{Theorem}  
Let $\mathbf{v}_k^i, \; \mathbf{u}_{kl}^i,\; \mathbf{y}_{kl}^i,\; i=1,\ldots,N$ is a sequence generated by algorithm \ref{algADMM} with $\mu>0$. Define $\bar{\mathbf{v}}_k=\sum_{i=1}^N \mathbf{v}_k^i/N$ and the same way $\bar{\mathbf{u}}_{kl},\;\bar{\mathbf{y}}_{kl}$. Then we have
\begin{eqnarray}\label{for function conv}
&f_k(\bar{\mathbf{v}}_k)-f_k(\mathbf{v}_k^*) &\leq \frac{\mu}{2N}\sum_{l\in \mathcal{B}_k}\|\mathbf{u}_{kl}^0- \mathbf{u}_{kl}^*\|^2_2, \notag\\
&\|\bar{\mathbf{v}}_k[l]-\bar{\mathbf{u}}_{kl}\|_2 &\leq \frac{1}{N}\left( \frac{2}{\mu}\|\mathbf{y}_{kl}^*\|_2+\|\mathbf{u}_{kl}^0-\mathbf{u}_{kl}^*)\|_2\right).
\end{eqnarray}
\end{thm}
The proof is presented in Appendix.


\begin{rem}
Note, that non-ergodic convergence rate of the original ADMM is $O(1/\sqrt{N})$. 
We consider the ergodic ADMM sequence $(\bar{\mathbf{v}}_k,\; \bar{\mathbf{u}}_{kl},\; \bar{\mathbf{y}}_{kl})$ in Theorem \ref{Theorem}, which has a higher convergence rate as compared to the non-ergodic ADMM approach. In addition, a convergence comparison for ergodic and non-ergodic methods can be found in \citet{Yin} paper. 
\end{rem}

\subsection{Accelerated Alternating Direction Method of Multiplier}
Although the ergodic ADMM exhibits good convergence, we suggest to apply the accelerated variant of ADMM from \citet{Goldstein} to our problem. This method is also covered in \citet{TanakaAADMM}, \citet{continuousAADMM}, \citet{lin2020accelerated} and others. The approach inherits the idea of \ag{accelerated} method from \citet{Nesterovpaper} (see \cite{guminov2019accelerated} for accelerated alternating procedures), in which a first-order minimization scheme with global convergence $O(1/N^2)$ for the class of \ag{smooth target} functions was presented.
Nesterov's accelerated method is a modification of gradient descent, which is accelerated by the step of overrelaxation. 

In the work of \citet{Goldstein}, the \ag{authors} assume that both functions from general problem (\ref{Lan_gen_prob}) are strongly convex. In our case, $g(\mathbf{y})$ is a zero function, so not strongly convex. Therefore, in \ag{theory}, the method \ag{have been applied for considered problem} may not yield \ag{any acceleration}, as was shown in Theorem 2 of \citet{Goldstein}. However, with the proper choice of parameters, the method can be
accelerated. Thus, we get an effective solution to the PSSE problem. \ag{Note, that in practice we observed for both of these methods ADMM and Accelerated ADMM linear rate of convergence! This is much better than can be explained by the best known theoretical results for such methods, e.g. see \cite{ouyang2015accelerated,Lan,lin2020accelerated}.}

The practical aspects is one of the main reasons to use the ADMM approach and its accelerated variants instead of not augmented and not alternating gradient type procedures and their accelerated variants. \pars{Moreover, it is preferable to use the A-ADMM, since the practical performance of the accelerated method, as proven through the experiments in the end of the paper, is better.} 
\begin{algorithm}
\caption{Accelerated ADMM}
\label{accelerated}
\begin{algorithmic}[1]
\Require matrix $H$, system state vector $\mathbf{z}$, parameter $\rho$, tolerances $\varepsilon_{primal}=10^{-3}, \; \varepsilon_{dual}=10^{-4}$, initial points $\mathbf{u}_{kl}^0 = 0, \, \hat{\mathbf{u}}_{kl}^0 = 0, \, \mathbf{y}_{kl}^0 = 0,\, \hat{\mathbf{y}}_{kl}^0 = 0$ for each $k, l \in K$
\State $i \leftarrow 0$
\Repeat
\State \begin{eqnarray*}
       &\mathbf{v}_k^i &= \min_{\mathbf{v}_k \in K} f_k(\mathbf{v}_k) + \sum_{l \in B_k} \hat{\mathbf{y}}_{kl}^{i-1}\left(\mathbf{v}_k[l] - \hat{\mathbf{u}}^{i-1}_{kl}\right) \\
       &&+ \sum_{l \in B_k} \frac{\rho}{2} \left\|\mathbf{v}_k[l] - \hat{\mathbf{u}}_{kl}^{i-1}\right\|_2^2
      \end{eqnarray*}
\State $\mathbf{u}_{kl}^i = \frac{1}{2}\left(\mathbf{v}_k^i[l] + \mathbf{v}_l^i[k]\right) $ 
\State $\mathbf{y}_{kl}^i = \hat{\mathbf{y}}^{i-1}_{kl} + \rho ( \mathbf{v}_k^i[l] - \mathbf{u}_{kl}^i)$ 
\State $\alpha_i = \frac{1+\sqrt{1+\alpha_{i-1}^2}}{2}$; 
\State $\hat{\mathbf{u}}_{kl}^i = \mathbf{u}_{kl}^i + \frac{\alpha_{i-1} - 1}{\alpha_i}\left(\mathbf{u}_{kl}^i- \mathbf{u}_{kl}^{i-1}\right)$; 
\State $\hat{\mathbf{y}}^i_{kl} = \mathbf{y}_{kl}^i + \frac{\alpha_{i-1} - 1}{\alpha_i}\left(\mathbf{y}_{kl}^i - \mathbf{y}_{kl}^{i-1}\right)$, 
\Until {$\|\mathbf{r}^i\|_2^2 < \varepsilon_{primal}$ and $\|\mathbf{s}^i\|_2^2 < \varepsilon_{dual}$}
\State $i \leftarrow i+1$
\end{algorithmic}
\end{algorithm}

\subsection{Discussion AADMM and other accelerated approaches}

\pars{Previously we compared the theoretical rate of convergence of ADMM with its different accelerated variants.} Nevertheless, without additional assumptions no theoretical acceleration can be reached. It was also confirmed in a recent book by \cite{lin2020accelerated}. 

Let us discuss now other competing approaches. For example, if we consider pure Accelerated Augmented Lagrange Multiplier Method (not with Alternating Direction), we can really guarantee acceleration $O(1/N^2)$ \cite{lin2020accelerated}. \pars{Unfortunately, this is achieved at cost of loosing an important feature of splitting computations between the areas (i.e. we loose the possibility to perform computations in a distributed manner). }

Another type of distributed approach to \eqref{optim_prob_for_distr_methods} is to rewrite the constraints $\mathbf{v}_k [l] = \mathbf{v}_l [k]$ without introducing of auxiliary $\mathbf{u}$. In this case we have a system of linear constraints with the Laplacian matrix (positive semi-definite) of some graph $W\mathbf{v}=0$.\footnote{It's important to note, that in the described below approaches (references) matrix $W$ has another definition than we propose here. But, the only property that is really required from $W$ in mentioned below approaches is to be a Laplacian matrix. This matrix corresponds to the graph with nodes enumerated according to $\mathbf{v}_k[l]$ for all admissible $k$, $l$. We draw an edge between two nodes $\mathbf{v}_k[l]$ and $\mathbf{v}_l[k]$ if $\mathbf{v}_k[l] = \mathbf{v}_l[k]$. We can also generalize this definition in case when three and more areas have common bus. Note, that if the situation described in the previous sentence is not a case, the degree of each node (the number of neighbors) is $0$ or $1$ and $\chi = 1$.
In general, the quadratic Laplacian matrix of such graph has size equal to the number of nodes. Diagonal element equals to the number of neighbors (for a considered node) and at corresponding row we put $-1$ at the positions that correspond to the neighbors (of considered node).} Thus, we can apply optimal decentralized distributed algorithms based on dual oracle, e.g. see \cite{scaman2017optimal,gasnikov2020modern,uribe2020dual,dvinskikh2019decentralized,gorbunov2019optimal,hendrikx2019asynchronous,hendrikx2020optimal,arjevani2020ideal}. In particular, Algorithm 5 from \cite{uribe2020dual}. This algorithm has almost the same iteration complexity and splitting power as ADMM, but requires $\tilde{O}(\chi/N^2)$ communication steps (iterations), where $\chi = \lambda_{\max}(W)/\lambda_{\min}^{+}(W)$, $\lambda_{\min}^{+}(W)$ is a minimal positive eigenvalue of matrix $W$. 

Moreover, for our problem formulation we do not have  dual-friendly oracle in general. That is the problems in line 3 of Algorithms~\ref{algADMM},~\ref{accelerated} are not necessarily simple (have explicit solution). Since that we can apply optimal primal decentralized distributed approaches by \cite{li2018sharp,dvinskikh2019decentralized,gorbunov2019optimal,xu2019accelerated,li2020revisiting,ye2020multi,kovalev2020optimal,hendrikx2020dual}, that still require $\tilde{O}(\chi/N^2)$ communications steps (and $\tilde{O}(1/N^2)$ $\nabla f_k(\mathbf{v}_k)$ calculations per each node). The main difference with the dual approach is that we do not need to solve the auxiliary problems of type line 3 in Algorithms~\ref{algADMM},~\ref{accelerated}. In primal approaches we will have instead of line 3 the explicit formulas that use $\nabla f_k(\mathbf{v}_k)$. 

So up to a factor of $\chi$ it seems that the ADMM approach is dominated by the described above approaches. But in general $\chi\sim \lambda_{\min}^{+}(W)^{-1}$ can be too large. In ADMM (to say more accurately in the ergodic version of ADMM that we consider \cite{lin2020accelerated}) the guaranteed rate of of convergence $O(1/N)$ does not depend on the network topology.\footnote{Note, that for ADMM  it is also possible to obtain the rates of convergence like $\tilde{O}(\chi/N^2)$ \cite{lin2020accelerated}. \pars{But this is not suitable for our case}.} Therefore, we prefer ADMM-type approaches, since they are more robust to network topology and typically work much better in practice
than the approaches mentioned above that are mostly well described by the developed theory.

If we had centralized distributed architecture we could use not only acceleration (in total number of $\nabla f_k(\mathbf{v_k})$ calculations): $1/N \to 1/N^2$, but also variance reduction acceleration that could additionally improve the complexity by a  factor of $\sqrt{K}$ \cite{Lan}. 

In this paper we do not assume that $\lambda_{\min}^{+}(W)$ is not small and(or) we have centralized architecture. That is why here we concentrate on ADMM-type methods, that have good practical relevance, e.g. see \cite{Boyd}. 

\pars{Note, that we are not the first who proposed to apply ADMM to the problem considered in the paper, see \cite{kekatos2013distributed}. Our contribution is in replacing the non-ergodic version of ADMM from \cite{kekatos2013distributed}, \cite{Kekatos} with the convergence $O(1/\sqrt{N})$ by the ergodic one (with rate of convergence $O(1/N)$), choosing the best accelerated version of ADMM (from practical considerations) and an attempt to explain that ADMM-type methods are not sensitive to network topology as other accelerated approaches.}

\section{Numerical Simulations}

The decentralized algorithms were tested on a 4,200-bus power grid synthetically
built from the IEEE 14- and 300-bus systems. Each of the 300
buses of the latter was assumed to be a different area, and was replaced by
a copy of the IEEE 14-bus grid. Additionally, every branch of the IEEE
300-bus grid was an inter-area line whose terminal buses are randomly selected
from the two incident to this line areas.

\begin{figure}[htb]
\begin{center}
\includegraphics[width=8.4cm]{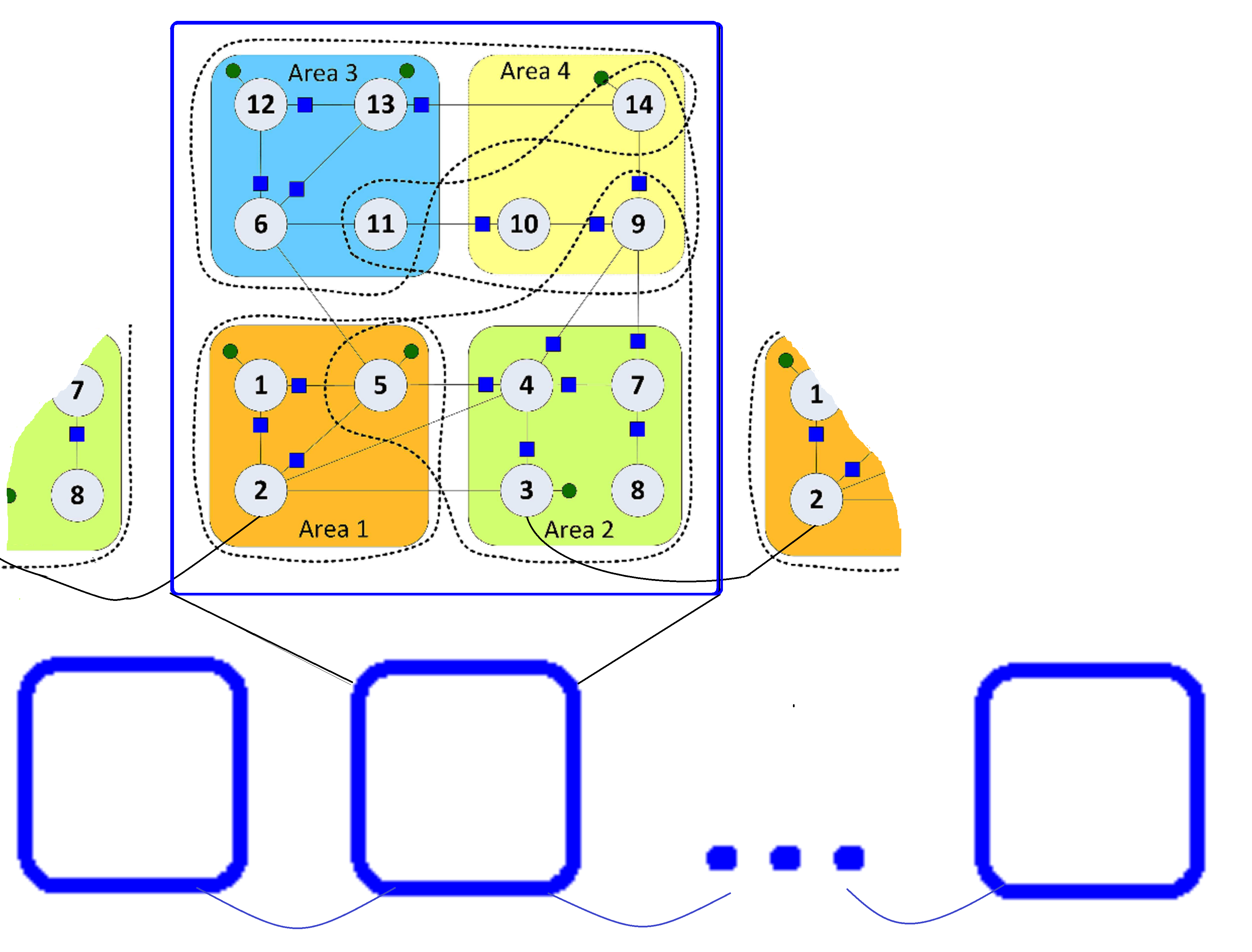}    
\caption{Hierarchical structure of the grid} 
\label{fig:GridEx}
\end{center}
\end{figure}

In this section, we validate and compare two distributed ADMM-based algorithms studied in the paper via simulations on power grids synthetically built from the IEEE 14-bus system. To provide a better comparison of the algorithms studied in the paper, we propose a scalable and flexible hierarchical topology of the network, organized as follows: the upper level of the network hierarchy is a chain (Example 2) or a ring (Examples 1,3) where each vertex was replaced by a copy of the IEEE 14-bus grid. Each $i$th 14-bus block is linked with its indexed neighboring blocks $i+1$ and $i-1$ by the lines as depicted in Fig.~\ref{fig:GridEx} (except for the first and last vertices in the case of chain topology of Example 2). The splitting of each block itself into four areas and the deployment of 22 sensors in each area was performed similarly to the example from Section~4, Fig.~2 in the work by \citet{Kekatos}.

In all examples we use the following parameters and convergence metrics: $\varepsilon_{dual} = 10^{-4},\; \varepsilon_{primal} = 10^{-3}$, where $\varepsilon_{dual}$ is a dual tolerance to achieve condition $\mathbf{u}_{kl}^{i} = \mathbf{u}_{kl}^{i+1}$ and $\varepsilon_{primal}$ is a primal tolerance to attain $\mathbf{v}_k[l] = \mathbf{u}_{kl}$ in each iteration. Convergence on direct variables achieved when all the state of areas reached the primal accuracy and convergence on dual variables checks the rapprochement for the consensus states of the network. The simulations were performed in PYTHON 3 using the CVXPY tools from \citet{CVX}.


\textbf{Example 1.:} 40 areas = 10 blocks linked into a ring, $\mu = 1,\; \rho = 2$.


Note, that in our topology the 14-bus blocks are linked quite weakly: each block can directly communicate with only two nearest neighbors, the cooperation with the other blocks is performed only indirectly. In the next example we decrease the connectivity by removing one line between the blocks and thus replacing a ring with a chain. This shows the importance of grid topology and its influence on the performance of distributed methods.


\begin{figure}[htb]
\begin{center}
\includegraphics[width=8.4cm]{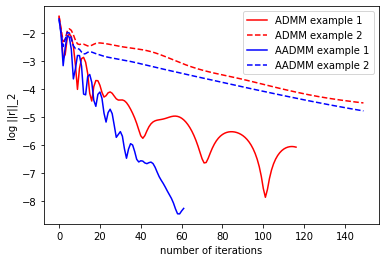}    
\includegraphics[width=8.4cm]{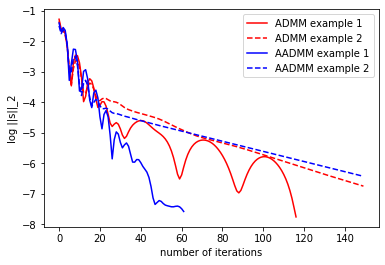}    
\caption{Primal and dual residuals versus number of iterations for Examples 1 and 2} 
\label{fig:Ex1}
\end{center}
\end{figure}

For each method, we show in Figure \ref{fig:Ex1} the logarithmic norm of primal residual for primal feasibility, which is defined as $\mathbf{r}_k^i = \mathbf{v}_k[l]^i - \mathbf{u}_{kl}^i$. And we show the logarithmic norm of the dual residual for the dual feasibility condition, which is defined as $\mathbf{s}_k^i = \mathbf{u}_{kl}^i - \mathbf{u}_{kl}^{i-1}$ in Figure \ref{fig:Ex1}. \\
We remark that the accelerated ADMM method significantly improves the convergence rate for primal and dual residuals in contrast with standard ADMM. In no more than 60 iterations, the accelerated method achieves convergence within the $\varepsilon_{primal} = 10^{-3}$, which is not accessible by the ADMM method after 100 iterations.

\begin{figure}[htb]
\begin{center}
\includegraphics[width=8.4cm]{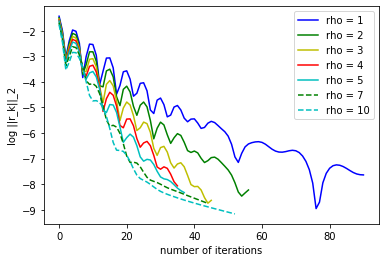}    
\includegraphics[width=8.4cm]{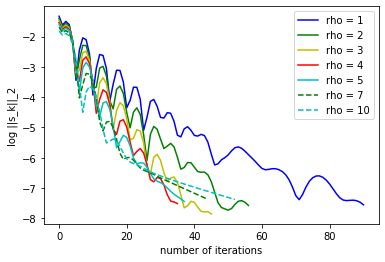}    
\caption{Primal and Dual residual versus number of iterations for Example 3} 
\label{fig:Ex2}
\end{center}
\end{figure}

Fig. \ref{fig:Ex2} illustrates sensitivity of the accelerated ADMM convergence to parameter $\rho$.

In the second example we demonstrate how significantly the changed topology may effect the convergence of the ADMM-based distributed methods.

\textbf{Example 2.:} 40 areas = 10 blocks linked into a chain, $\mu = 1,\; \rho = 2$.

As seen in Fig. 1, the convergence of both methods sharply worsened due to a break in communication. Which clearly shows the importance of network topology.

\textbf{Example 3.:} performance comparison within ring topology for a different number of areas\\
In Table \ref{comparison} one can see the comparative performance of ADMM and A-ADMM. For each method and each number of areas, we selected the best parameters $\mu$ and $\rho$ in terms of lower number of iterations to satisfy the conditions of stopping criterion.
\begin{table}[H]
\caption{Comparison of ADMM and A-ADMM}\label{comparison}
\begin{center}
\begin{tabular}{ccccc}
area & $\mu$ & ADMM, & $\rho$ & A-ADMM,\\
number & & iteration & & iteration\\
& & number && number\\
\hline
4 & 8 & 18 & 10 & 14 \\
20 & 5 & 23 & 6 & 20 \\
40 & 4 & 43 & 4 & 35 \\
60 & 2 & 65 & 3 & 44 \\
80 & 4 & 84 & 4 & 56 \\
100 & 2 & 115 & 2 & 82 \\
120 & 4 & 109 & 4 & 80 \\
\hline
\end{tabular}
\end{center}
\end{table} 

Note, that our convergence criterion comprises two components related to convergence of 1) the shared variables (consensus term) and 2) the local states of the areas. Due to this fact, it is always a trade-off how to set parameters $\mu$ and $\rho$. The table shows that in small-scale networks the best parameter values are about 8-10. On the other hand, if the network is large and sparse, we set smaller weights for the consensus terms, thus adjusting not only the common states of the interacting areas, but also the states within each area. 


%
%

\section{Conclusion}


In our paper, we analyze the performance of ADMM-based distributed state estimation approaches for power systems with linear measurement models and also discuss the applicability of some novel accelerated optimization methods. Our contribution is threefold:
\begin{enumerate}
    \item improvement of the conventional non-ergodic ADMM for distributed PSSE, guaranteeing the rate of convergence  $O(1/N)$ (ergodic) instead $O(1/\sqrt{N})$ (non-ergodic);
    \item reformulation of an accelerated version of ADMM in a distributed fashion, analysis of its theoretical and practical performance and comparison of ADMM and A-ADMM;
    \item discussion of other possible ways to speed up the convergence of distributed SE in cases where extra information about communication topology is available. 
\end{enumerate}

The future plans are to apply the methods discussed in the paper to 1) the case of asynchronous communications and 2) the nonlinear model of AC power flows.






\section{Appendix} 
We proof the Theorem 1 in a similar way as in  \citet{Lan}.
Let the primal-dual gap function $Q$ be
\begin{eqnarray}
\label{gap function}
&Q(\mathbf{t}_k^i, \mathbf{t}_k)&= f_k(\mathbf{v}_k^i)+\sum_{l \in \mathcal{B}_k} \langle \mathbf{y}_{kl}, \mathbf{v}_k^{i}[l]-\mathbf{u}_{kl}^i\rangle -f_k(\mathbf{v}_k)\notag \\
&&-\sum_{l \in \mathcal{B}_k} \langle \mathbf{y}_{kl}^{i}, \mathbf{v}_k[l]-\mathbf{u}_{kl}\rangle,
\end{eqnarray}
where $\mathbf{t}_k^i \equiv (\mathbf{v}_k^i,\mathbf{u}_{kl}^i, \mathbf{y}_{kl}^i )$ and $\mathbf{t}_k \equiv (\mathbf{v}_k,\mathbf{u}_{kl}, \mathbf{y}_{kl})$. 

\begin{pf}    

1) Due to the optimality conditions of the 
\begin{align*}
\mathbf{v}_k^i &= \min_{\mathbf{v}_k \in K} f(\mathbf{v}_k) + \sum_{l \in B_k} \mathbf{y}_{kl}^{i-1}(\mathbf{v}_k[l] - \mathbf{u}^{i-1}_{kl}) \notag \\
&+ \sum_{l \in B_k} \frac{\mu}{2} \|\mathbf{v}_k[l] - \mathbf{u}_{kl}^{i-1}\|_2^2,
\end{align*}
we obtain \begin{eqnarray}
&&f_k(\mathbf{v}_k^i)-f_k(\mathbf{v}_k) \notag \\
&+&\sum_{l\in \mathcal{B}_k}\langle \mathbf{y}_{kl}^{i-1},\mathbf{v}_k^i[l]-\mathbf{u}_{kl}^{i-1}\rangle -\sum_{l\in \mathcal{B}_k} \langle \mathbf{y}_{kl}^{i-1},\mathbf{v}_k[l]-\mathbf{u}_{kl}^{i-1}\rangle \notag\\
&+& \frac{\mu}{2}\sum_{l\in \mathcal{B}_k}\|\mathbf{v}_k^i[l]-\mathbf{u}_{kl}^{i-1}\|_2^2-\frac{\mu}{2}\sum_{l\in \mathcal{B}_k}\|\mathbf{v}_k[l]-\mathbf{u}_{kl}^{i-1}\|_2^2 \notag \\
&\geq& f_k(\mathbf{v}_k^i)-f_k(\mathbf{v}_k)+\sum_{l\in\mathcal{B}_k}\langle \mathbf{y}_{kl}^{i-1},\mathbf{v}_k^i[l]-\mathbf{v}_k[l]\rangle \notag\\
&+& \frac{\mu}{2}\sum_{l\in \mathcal{B}_k} \langle 2(\mathbf{v}_k^i[l]-\mathbf{u}_{kl}^{i-1}),\mathbf{v}_k^i[l]-\mathbf{v}_k[l]\rangle  = \notag \\
&=&f_k(\mathbf{v}_k^i)-f_k(\mathbf{v}_k)\notag \\
&+&\sum_{l\in \mathcal{B}_k}\langle \mathbf{y}_{kl}^{i-1}+\mu(\mathbf{v}_k^i[l]-\mathbf{u}_{kl}^{i-1}),\mathbf{v}_k^i[l]-\mathbf{v}_k[l]\rangle \leq 0,
\end{eqnarray}
which together with the equation $\mathbf{y}_{kl}^i = \mathbf{y}^{i-1}_{kl} + \mu ( \mathbf{v}_k^i[l] - \mathbf{u}_{kl}^i)$ represents that
\[
f_k(\mathbf{v}_k^i)-f_k(\mathbf{v}_k)\leq - \sum_{l\in \mathcal{B}_k} \langle \mathbf{y}_{kl}^i+\mu (\mathbf{u}_{kl}^i-\mathbf{u}_{kl}^{i-1}), \mathbf{v}_k^i[l]-\mathbf{v}_k[l]\rangle.
\]
2) Take these obtained relationships and definition of the duality gap (\ref{gap function}), we get
\begin{eqnarray}
&&Q(\mathbf{t}_k^i,\mathbf{t}_k)\leq  -\sum_{l\in \mathcal{B}_k}\langle \mathbf{y}_{kl}^i+\mu (\mathbf{u}_{kl}^i-\mathbf{u}_{kl}^{i-1}),\mathbf{v}_k^i[l]-\mathbf{v}_k[l]\rangle \notag \\
&-&\sum_{l\in \mathcal{B}_k}\langle \mathbf{y}_{kl}^i,\mathbf{v}_k[l]-\mathbf{v}_{kl}\rangle +\sum_{l\in \mathcal{B}_k}\langle \mathbf{y}_{kl}, \mathbf{v}_k^i[l]-\mathbf{u}_{kl}^i \rangle \notag\\
&=&\sum_{l\in \mathcal{B}_k}\langle \mathbf{y}_{kl}^i,\mathbf{u}_{kl}-\mathbf{v}_k^i[l] \rangle + \sum_{l\in \mathcal{B}_k}\mu\langle \mathbf{u}_{kl}^{i-1}-\mathbf{u}_{kl}^i,\mathbf{v}_k^i[l]-\mathbf{v}_k[l]\rangle \notag\\
&+&\sum_{l\in \mathcal{B}_k}\langle \mathbf{y}_{kl}, \mathbf{v}_k^i[l]-\mathbf{u}_{kl}^i\rangle 
=\sum_{l\in \mathcal{B}_k}\langle \mathbf{y}_{kl}-\mathbf{y}_{kl}^i, \frac{1}{\mu}(\mathbf{y}_{kl}^i-\mathbf{y}_{kl}^{i-1})\rangle \notag \\
&+& 
\sum_{l\in \mathcal{B}_k} \mu \langle \mathbf{u}_{kl}^{i-1}-\mathbf{u}_{kl}^i, \mathbf{v}_k^i[l]-\mathbf{v}_k[l]\rangle.
\end{eqnarray}
3) Note, that
\begin{eqnarray*}
&&2\langle \mathbf{y}_{kl}-\mathbf{y}_{kl}^i, \mathbf{y}_{kl}^i-\mathbf{y}_{kl}^{i-1} \rangle \\
&=&\|\mathbf{y}_{kl}-\mathbf{y}_{kl}^{i-1}\|^2_2-\|\mathbf{y}_{kl}-\mathbf{y}_{kl}^i\|^2_2-\|\mathbf{y}_{kl}^{i-1}-\mathbf{y}_{kl}^i\|^2_2; \\
&&2\langle \mathbf{u}_{kl}^{i-1}-\mathbf{u}_{kl}^i,\mathbf{v}_k^i[l]-\mathbf{v}_k[l]\rangle \\
&=&\|\mathbf{v}_k[l]-\mathbf{u}_{kl}^{i-1}\|^2_2-\|\mathbf{v}_k[l]-\mathbf{u}_{kl}^i\|^2_2+\|\mathbf{v}_k^i[l]-\mathbf{u}_{kl}^i\|^2_2 \\
&-&\|\mathbf{v}_k^i[l]-\mathbf{u}_{kl}^{i-1}\|^2_2=\|\mathbf{v}_k[l]-\mathbf{u}_{kl}^{i-1}\|^2_2-\|\mathbf{v}_k[l]-\mathbf{u}_{kl}^i\|^2_2 \\
&+&\frac{1}{\mu^2}\|\mathbf{y}_{kl}^i-\mathbf{y}_{kl}^{i-1}\|^2_2-\|\mathbf{v}_k^i[l]-\mathbf{u}_{kl}^{i-1}\|^2_2. \\
\end{eqnarray*}

As a result, we have
\begin{eqnarray}
&&Q(\mathbf{t}_k^i,\mathbf{t}_k)\leq \frac{1}{2\mu}\sum_{l\in \mathcal{B}_k}\left( \|\mathbf{y}_{kl}-\mathbf{y}_{kl}^{i-1}\|^2_2-\|\mathbf{y}_{kl}-\mathbf{y}_{kl}^i\|^2_2\right) \notag \\
&&-\frac{1}{2\mu}\sum_{l\in \mathcal{B}_k}\left(\|\mathbf{y}_{kl}^{i-1}-\mathbf{y}_{kl}^i\|^2_2 \right) \notag\\
&&+\frac{\mu}{2}\sum_{l\in \mathcal{B}_k}\left( \|\mathbf{v}_k[l]-\mathbf{u}_{kl}^{i-1}\|^2_2-\|\mathbf{v}_k[l]-\mathbf{u}_{kl}^i\|^2_2\right) \notag\\
&&+\frac{\mu}{2}\sum_{l\in \mathcal{B}_k}\left(\frac{1}{\mu^2}\|\mathbf{y}_{kl}^i-\mathbf{y}_{kl}^{i-1}\|^2_2-\|\mathbf{v}_k^i[l]-\mathbf{u}_{kl}^{i-1}\|^2_2\right) \notag\\
&=& \frac{1}{2\mu}\sum_{l\in \mathcal{B}_k}\left(\|\mathbf{y}_{kl}-\mathbf{y}_{kl}^{i-1}\|^2_2-\|\mathbf{y}_{kl}-\mathbf{y}_{kl}^i\|^2_2\right) \notag\\
&&+\frac{\mu}{2}\sum_{l\in \mathcal{B}_k} ( \|\mathbf{v}_k[l]-\mathbf{u}_{kl}^{i-1}\|^2_2-\|\mathbf{v}_k[l]-\mathbf{u}_{kl}^i\|^2_2 \notag\\
&&-\|\mathbf{v}_k^i[l]-\mathbf{u}_{kl}^{i-1}\|^2_2 ). \notag
\end{eqnarray}
4) Obtained inequality summing up for all $i=1,\ldots,N$, we have
\begin{eqnarray}\label{sum Q by i}
&&\sum_{i=1}^N Q(\mathbf{t}_k^i,\mathbf{t}_k)\leq \frac{1}{2\mu}\left(\|\mathbf{y}_{kl}^{0}-\mathbf{y}_{kl}\|^2_2-\|\mathbf{y}_{kl}^N-\mathbf{y}_{kl}\|^2_2\right) \notag\\
&+&\frac{\mu}{2}\left(\|\mathbf{u}_{kl}^0-\mathbf{u}_{kl}\|^2_2-\|\mathbf{u}_{kl}^N-\mathbf{u}_{kl}\|^2_2\right) \geq 0.
\end{eqnarray}
5) Substitute $ \mathbf{t}_k$ by $ \mathbf{t}_k^* = (\mathbf{v}_k^*, \mathbf{u}_{kl}^*, \mathbf{y}_{kl}^*) $ in (\ref{sum Q by i}) and take into account, that $Q (\bar{\mathbf{t}}_k, \mathbf{t}^*)\geq 0 $, we see
\[
\|\mathbf{y}_{kl}^N-\mathbf{y}_{kl}^*\|^2_2\leq \|\mathbf{y}_{kl}^{0}-\mathbf{y}_{kl}^*\|^2_2+\mu^2\|\mathbf{u}_{kl}^0-\mathbf{u}_{kl}^*\|^2_2.
\]
Therefore
\begin{eqnarray}\label{in theor lambda}
&&\|\mathbf{y}_{kl}^{N}-\mathbf{y}_{kl}^0\|_2\leq \|\mathbf{y}_{kl}^{0}-\mathbf{y}_{kl}^*\|^2_2+\|\mathbf{y}_{kl}^{N}-\mathbf{y}_{kl}^*\|_2 \notag\\
&&\leq 2\|\mathbf{y}_{kl}^{0}-\mathbf{y}_{kl}^*\|_2+\mu\|\mathbf{v}_{kl}^0-\mathbf{u}_{kl}^*\|_2.
\end{eqnarray}
6) Let $\mathbf{t}_k=(\mathbf{v}_k^*,\mathbf{u}_{kl}^*,\mathbf{y}_{kl})$ in (\ref{sum Q by i}), then note that
\begin{eqnarray*}
&&\frac{1}{N}\sum_{i=1}^N Q(\mathbf{t}_k^i,(\mathbf{v}_k^*,\mathbf{u}_{kl}^*,\mathbf{y}_{kl}))\geq Q(\bar{\mathbf{t}}_k,(\mathbf{v}_k^*,\mathbf{u}_{kl}^*,\mathbf{y}_{kl})) \\
&=&f_k(\bar{\mathbf{v}}_k)-f_k(\mathbf{v}_k^*)+\sum_{l\in \mathcal{B}_k}\langle \mathbf{y}_{kl},\bar{\mathbf{v}}_k[l]-\bar{\mathbf{u}}_{kl}\rangle,
\end{eqnarray*}
and

\begin{eqnarray*}
&\frac{1}{2}(||\mathbf{y}_{kl}^{0}-\mathbf{y}_{kl}||^2_2-||\mathbf{y}_{kl}^{N}-\mathbf{y}_{kl}||^2_2)&=
\frac{1}{2}(||\mathbf{y}_{kl}^0||^2_2-||\mathbf{y}_{kl}^N||^2_2)\\
&&-\langle \mathbf{y}_{kl},\mathbf{y}_{kl}^0-\mathbf{y}_{kl}^N\rangle.
\end{eqnarray*}

Hence, we have
\begin{eqnarray*}
&& f_k(\bar{\mathbf{v}}_k)-f_k(\mathbf{v}_k^*)+\sum_{l\in \mathcal{B}_k}\langle \mathbf{y}_{kl},\bar{\mathbf{v}}_k[l]-\bar{\mathbf{u}}_{kl} \rangle \\
&\leq& \frac{1}{N}\sum_{l\in \mathcal{B}_k}\frac{1}{\mu}\Bigl[\frac{1}{2}(\|\mathbf{y}_{kl}^0\|^2_2-\|\mathbf{y}_{kl}^N\|^2_2) - \langle \mathbf{y}_{kl}, \mathbf{y}_{kl}^0-\mathbf{y}_{kl}^N\rangle \Bigr] \\
&+&\frac{\mu}{2N}\sum_{l\in \mathcal{B}_k}\Bigl[\|\mathbf{u}_{kl}^0-\mathbf{u}_{kl}^*\|^2_2-\|\mathbf{u}_{kl}^N-\mathbf{u}_{kl}^*\|^2_2\Bigr],
\end{eqnarray*}    
7) Therefore,
\begin{eqnarray*}
&&f_k(\bar{\mathbf{v}}_k)-f_k(\mathbf{v}_k^*)+\sum_{l\in \mathcal{B}_k}\langle \mathbf{y}_{kl},\bar{\mathbf{v}}_k[l]-\bar{\mathbf{u}}_{kl} + \frac{1}{\mu N}(\mathbf{y}_{kl}^0-\mathbf{y}_{kl}^N) \rangle \\
&&\leq\frac{1}{2\mu N}\sum_{l\in \mathcal{B}_k}\Bigl( \|\mathbf{y}_{kl}^0\|^2_2-\|\mathbf{y}_{kl}^N\|^2_2\Bigr)   \\
&&+\frac{\mu}{2N}\sum_{l\in \mathcal{B}_k}\Bigl( \|\mathbf{u}_{kl}^0-\mathbf{u}_{kl}^*\|^2_2-\|\mathbf{u}_{kl}^N-\mathbf{u}_{kl}^*\|^2_2 \Bigr)\\
&&\leq \frac{1}{2N}\Bigl( \frac{1}{\mu}\sum_{l\in \mathcal{B}_k}\|\mathbf{y}_{kl}^0\|_2^2+\mu \sum_{l\in \mathcal{B}_k}\|\mathbf{u}_{kl}^0-\mathbf{u}_{kl}\|_2^2 \Bigr)\\
&&= \frac{\mu}{2N}\sum_{l\in \mathcal{B}_k}\|\mathbf{u}_{kl}^0-\mathbf{u}_{kl}\|_2^2,
\end{eqnarray*}
for all $\mathbf{u}_{kl}, l \in\mathcal{B}_k, k=1, \ldots, K$. In the last equality, we took into account that we initialized $\mathbf{y}_{kl}^0=0$. Further, from this relation follows, that $\bar{\mathbf{v}}_k[l]-\bar{\mathbf{v}}_{kl}+\frac{1}{\mu N}(\mathbf{y}_{kl}^0-\mathbf{y}_{kl}^k)=0$,  and hence the first equation in (\ref{for function conv}) is satisfied. Also, from the last and (\ref{in theor lambda}) follow the second equation from (\ref{for function conv}).  
\begin{eqnarray*}
&\|\bar{\mathbf{v}}_k[l]-\bar{\mathbf{u}}_{kl}\|_2 &\leq \frac{1}{\mu N} \Bigl[2\|\mathbf{y}_{kl}^0-y_{kl}^*\|_2+\mu\|\mathbf{u}_{kl}^0-\mathbf{u}_{kl}^*\|_2\Bigr]  \\
&&=\frac{1}{N}\Bigl[ \frac{2}{\mu} \|\mathbf{y}_{kl}^0-\mathbf{y}_{kl}^*\|_2+\|\mathbf{u}_{kl}^0-\mathbf{u}_{kl}^*\|_2\Bigr]\\
&&=\frac{1}{N}\Bigl[ \frac{2}{\mu} \|\mathbf{y}_{kl}^*\|_2+\|\mathbf{u}_{kl}^0-\mathbf{u}_{kl}^*\|_2\Bigr].
\end{eqnarray*}

\end{pf}

\end{document}